\documentstyle[12pt]{article}
\title{Necessary And Sufficient Conditions For Existence of the LU Factorization  of an Arbitrary Matrix.}
\author{Adviser: Charles R. Johnson \\ Department of
Mathematics\\ College of William and Mary\\Williamsburg, VA \and Pavel
Okunev\thanks{This research was partially supported by NSF REU grant
DMS-96-19577 and by the College of William and Mary.} \\  Department of
Mathematics\\
University of Connecticut\\ Storrs, CT}
\date{July 31, 1997}
\begin{document}
\maketitle

\begin{abstract}
If $A$ is an n-by-n matrix over a field $F$ ($A\in M_{n}(F)$), then $A$ is
said to ``have an LU factorization'' if there exists a lower
triangular matrix $L\in M_{n}(F)$ and an upper triangular matrix $U\in
M_{n}(F)$ such that $$A=LU.$$ We give necessary and sufficient
conditions for LU factorability of a matrix. Also simple algorithm for
computing an LU factorization is given. It is an extension of the Gaussian elimination algorithm to the
case of not necessarily invertible matrices. We consider possibilities to
factors a matrix that does not have an LU factorization as the product
of an ``almost lower triangular'' matrix and an ``almost upper
triangular'' matrix. There are many ways to formalize what almost
means. We consider some of them and derive necessary and sufficient
conditions. Also simple algorithms for  computing of  an ``almost LU
factorization'' are given.

\end{abstract}

\newtheorem{theorem}{Theorem}
\newtheorem{lemma}{Lemma}
\newtheorem{definition}{Definition}
\newtheorem{corollary}{Corollary}
\newtheorem{note}{Note}
If $A$ is an n-by-n matrix over a field $F$ ($A\in M_{n}(F)$), then $A$ is
said to ``have an LU factorization'' if there exists a lower
triangular matrix $L\in M_{n}(F)$ and an upper triangular matrix $U\in
M_{n}(F)$ such that $$A=LU.$$ Historically, such factorizations have
been extensively used in the solution of the linear system $Ax=b$,
especially when a series of right hand sides $b$ is presented. For, if
$A=LU$, then each solution $x$ is a solution to $$Ux=y$$ for some
solution $y$ to $$Ly=b,$$ and if $L$ is invertible, $y$ is uniquely
determined; and if $U$ is invertible, any solution $y$ to $Ly=b$ will
give a solution $x$ to $Ax=b$ via $Ux=y$. Thus, much of the discussion
of LU factorization has occurred in the computational literature .
However, LU factorization has also emerged as a useful theoretical
tool.

Though an LU factorization is quite useful when it exists,
unfortunately not every matrix has one. The simplest example is
$$A=\left[\begin{array}{cc}
0&1\\
1&0
\end{array}\right].$$
If $A=LU=\left[\begin{array}{cc}
l_{11}&0\\
l_{21}&l_{22}
\end{array}\right]\left[\begin{array}{cc}
u_{11}&u_{12}\\
0&u_{22}
\end{array}\right]$, then $L$ and $U$ would have to be invertible
because $A$ is
(i.e. $l_{11},u_{11},l_{22},u_{22}\neq 0$), but this would contradict
$l_{11}u_{11}=a_{11}=0$.

Which $A\in M_{n}(F)$, then, do have an LU factorization? We first
observe a necessary condition on $A$ by writing a supposed
factorization in partitioned form. Suppose
$$A=\left[\begin{array}{cc}
A_{11}&A_{12}\\
A_{21}&A_{22}
\end{array}\right]$$ with $A_{11}$ k-by-k, and that
$$A=\left[\begin{array}{cc}
L_{11}&0\\
L_{21}&L_{22}
\end{array}\right]\left[\begin{array}{cc}
U_{11}&U_{12}\\
0&U_{22}
\end{array}\right]$$ partitioned conformally. We then have
$$[A_{11}A_{12}]=L_{11}[U_{11}U_{12}]$$ and $$\left[\begin{array}{cc}
A_{11}\\
A_{21}
\end{array}\right]=\left[\begin{array}{cc}
L_{11}\\
L_{21}
\end{array}\right]U_{11}$$ and $$A_{11}=L_{11}U_{11}$$ We use well
known inequalities involving ranks of products. If $Z=XY$ with $X$
having $k$ columns and $Y$ having $k$ rows,
$$rankX+rankY-k\leq rankZ\leq min\{ rankX,rankY\}$$
First applying the right hand inequality we have $$rank[A_{11}A_{12}]\leq
rank L_{11}$$ and $$rank\left[\begin{array}{cc}
A_{11}\\
A_{21}
\end{array}\right]\leq rankU_{11}$$
Then applying the left hand inequality we get
$$rankL_{11}+rankU_{11}-k\leq rankA_{11}$$
Combining, we obtain $$rank[A_{11}A_{12}]+rank\left[\begin{array}{cc}
A_{11}\\
A_{21}
\end{array}\right]-k\leq rankA_{11}$$ which we may write as
$$rankA_{11}+k\geq rank[A_{11}A_{12}]+rank\left[\begin{array}{cc}
A_{11}\\
A_{21}
\end{array}\right]$$ this must hold for each $k=1,\ldots ,n-1$. We can
also write these conditions as

\begin{eqnarray}
\label{main}
rankA[\{1 \ldots k\}]+k\geq rankA[\{1 \ldots k\},\{1 \ldots n\}]+rankA[\{1
\ldots n\},\{1 \ldots k\}]
\end{eqnarray}
for all $k=1,\ldots ,n$.

Our main result is that these conditions are also sufficient.

\begin{theorem}

\label{lu}

The matrix $A=(a_{ij})\in M_{n}(F)$ has an LU factorization iff it satisfies
conditions (\ref{main}).

\end{theorem}

{\em Proof.} Necessity was proved above. We will prove sufficiency by
induction on n. In the
case $n=1$ the theorem is trivially true because any 1-by-1 matrix has
an LU factorization. In particular, if $A=[0]$ we let $L=[0]$ and
$U=[0]$. Though there are many other LU factorizations of the zero matrix
we choose this particular one for the purposes of theorem \ref{long}.

Suppose that conditions (\ref{main}) are sufficient for $n\leq m$. We
want to prove that conditions (\ref{main}) are also sufficient for $n=m+1$.

The idea of the proof is as follows. We want to find a
lower triangular matrix L and an upper triangular matrix U such that
$A=LU$. Suppose $$A=\left[ \begin{array}{cc}
a_{11}&A_{12}\\
A_{21}&A_{22}
\end{array}\right]$$ in which $a_{11}$ is the (1,1) entry of A, and
that $$A=LU=\left[\begin{array}{cc}
l_{11}&0\\
L_{21}&L_{22}
\end{array}\right]\left[ \begin{array}{cc}
u_{11}&U_{12}\\
0&U_{22}
\end{array}\right]$$ partitioned conformally.

Block multiplication gives us
\begin{eqnarray}
a_{11}&=&l_{11}u_{11} \label{eq:block1}\\
A_{21}&=&L_{21}u_{11}\\
A_{12}&=&l_{11}U_{12} \label{eq:block3}\\
A_{22}&=&L_{21}U_{12}+L_{22}U_{22}. \label{eq:block4}
\end{eqnarray}

We will prove that it is possible to choose the first column of L and the
first row of U in
such a way that equations \ref{eq:block1}--\ref{eq:block3} are
satisfied and additionally that the (n-1)-by-(n-1) matrix
$$B=A_{22}-L_{21}U_{12}$$ satisfies the conditions (\ref{main}). We
call $B$ ``a factor complement'' of $A$.
Then by the inductive hypothesis matrix B has an LU factorization or
$B=L_{0}U_{0}$ with $L_{0}$ lower triangular and $U_{0}$ upper
triangular. We let $$L_{22}=L_{0}$$ and $$U_{22}=U_{0}.$$ Notice that
a matrix L defined this way is lower triangular, and a matrix U defined
this way is upper triangular. Also equations
\ref{eq:block1}--\ref{eq:block4} are satisfied which implies that
$A=LU$.

We will now describe how to choose the first column of $L$ and the
first row of $U$. We need to consider several cases.

{\em Case 1: $a_{11}\neq 0$.} In this case we let
$$L[\{1 \ldots n\},\{1\}]=A[\{1 \ldots n\},\{1\}]$$ and $$U[\{1\},\{1
\ldots n\}]=a_{11}^{-1}A[\{1\},\{1 \ldots n\}]$$ which gives us
\begin{eqnarray}
l_{11}&=&a_{11} \nonumber \\
L_{21}&=&A_{21} \nonumber \\
u_{11}&=&1 \nonumber \\
U_{12}&=&a_{11}^{-1}A_{12}. \nonumber
\end{eqnarray}
Equations (\ref{eq:block1})--(\ref{eq:block3}) are obviously satisfied by
this choice of the first column of $L$ and the first row of $U$. It
remains to prove that the matrix $B=A_{22}-L_{21}U_{12}$ satisfies
conditions (\ref{main}).

Let us consider the matrix
\begin{equation}
\label{nzero}
C=A-\left[\begin{array}{c}
a_{11}\\
L_{21}
\end{array}\right] [0 \ U_{12}]=\left[\begin{array}{cc}a_{11}&0\\
A_{21}&A_{22}-L_{21}U_{12}
\end{array}\right]=\left[\begin{array}{cc}
a_{11}&0\\
A_{21}&B
\end{array}\right].
\end{equation}
Matrix $C$ satisfies conditions (\ref{main}). Indeed, let us fix k
($1\leq k \leq n$). Then matrix $C[\{1 \ldots k \}]$ can be obtained from
matrix $A[\{1 \ldots k \}]$ by applying several type-3 elementary column
operations. Therefore $$rankC[\{1 \ldots  k\}]=rankA[\{1 \ldots
 k \}].$$ In particular, $$rankC=rankA$$ Similarly, $$rankC[\{1 \ldots k\},\{1
\ldots n\}]=rankA[\{1 \ldots
k\},\{1 \ldots n\}]$$ and $$rankC[\{1 \ldots n\},\{1 \ldots k\}]=rankA[\{1
\ldots
n\},\{1 \ldots k\}].$$ Since $A$ satisfies conditions (\ref{main}) we
conclude that $C$ also satisfies conditions (\ref{main}).

From the formula (\ref{nzero}) it now follows that matrix $B$
satisfies conditions (\ref{main}). Indeed, let us fix k ($1 \leq k
\leq n-1$). Because $c_{11}=a_{11}\neq 0$ and $c_{12}=\cdots=c_{1,k+1}=0$
we have $$rankC[\{1\ldots n\},\{2 \ldots k+1\}]=rankC[\{1\ldots n\}, \{1 \ldots
k+1\}]-1$$ Also because $c_{12}=\cdots=c_{1,k+1}=0$ we see that
$$rankC[\{ 1 \ldots n\},\{ 2 \ldots k+1\}]]=rankC[\{ 2 \ldots n\},\{ 2 \ldots
k+1\}],$$
but $C[\{ 2 \ldots n\},\{ 2 \ldots k+1\}]=B[\{1\ldots n-1\},\{1 \ldots k\}]$
and therefore
\begin{eqnarray}
rankB[\{1 \ldots (n-1)\},\{1 \ldots k\}]=rankC[\{1 \ldots
n\},\{1 \ldots (k+1)\}]-1 \label{eq:nz1}
\end{eqnarray}

Similarly,
\begin{eqnarray}
&rankB[\{1 \ldots k\},\{1 \ldots (n-1)\}]=rankC[\{1 \ldots (k+1)\},\{1
\ldots n\}]-1 \\
&rankB[\{1 \ldots  k\}]=rankC[\{1 \ldots (k+1) \}]-1 \label{eq:nz3}
\end{eqnarray}
In particular, $$rankB=rankC-1=rankA-1$$
Also since $C$ satisfies conditions (\ref{main}) we have
\begin{eqnarray}
rankC[\{1 \ldots (k+1) \}]+(k+1)&\geq&rankC[\{1 \ldots (k+1)\},\{1
\ldots n\}]+ \nonumber \\
 & &rankC[\{1 \ldots n\},\{1 \ldots (k+1)\}]. \label{eq:nz4}
\end{eqnarray}
From (\ref{eq:nz1})--(\ref{eq:nz4}) we conclude that
\begin{eqnarray}
rankB[\{1 \ldots k \}]+k&\geq&rankB[\{1 \ldots k\},\{1 \ldots
(n-1)\}]+\nonumber \\
 & &rankB[\{1 \ldots (n-1)\},\{1 \ldots k\}]. \nonumber
\end{eqnarray} In other words $B$ satisfies conditions
(\ref{main}).

{\em Case 2: $a_{11}=0$.} If $a_{11}=0$, then conditions (\ref{main})
for $k=1$ imply that either the first row or the first column of $A$
or both are zero. We first consider the case when the first row of $A$
is equal to zero but the first column is not zero. Let $i$ be the
smallest integer such that $a_{1i}\neq 0$. We let
$$L[\{1 \ldots n\},\{1\}]=A[\{1 \ldots n\},\{1\}]$$ and $$U[\{1\},\{1
\ldots n\}]=a_{1i}^{-1}A[\{i\},\{1 \ldots n\}]$$ which gives us
\begin{eqnarray}
l_{11}&=&0 \nonumber \\
L_{21}&=&A_{21} \nonumber \\
u_{11}&=&1 \nonumber \\
U_{12}&=&a_{1i}^{-1}A[\{i\},\{2 \ldots n\}] \nonumber
\end{eqnarray}

Again equations (\ref{eq:block1})--(\ref{eq:block3}) are obviously satisfied by
this choice of the first column of $L$ and the first row of $U$. It
remains to prove that the matrix
$B=A_{22}-L_{21}U_{12}$ satisfies
conditions (\ref{main}).

Let us consider the matrix
\begin{equation}
\label{zero}
C=A-\left[\begin{array}{c}
0\\
L_{21}
\end{array}\right] [0 \ U_{12}]=\left[\begin{array}{cc}0&0\\
A_{21}&A_{22}-L_{21}U_{12}
\end{array}\right]=\left[\begin{array}{cc}
0&0\\
A_{21}&B
\end{array}\right].
\end{equation}
Exactly as in case 1 we prove that matrix $C$ satisfies
conditions (\ref{main}), and $$rankC=rankA$$

From the formula (\ref{zero}) it now follows that matrix $B$
satisfies conditions (\ref{main}). Indeed, let us fix k ($1 \leq k
\leq n-1$). Because $c_{i1}=a_{1i}\neq 0$ and $c_{i2}=\cdots=c_{i,k+1}=0$
we have $$rankC[\{1 \ldots n\},\{2 \ldots k+1\}]=rankC[\{1 \ldots
 n\},\{1 \ldots k+1\}]-1$$ Also because $c_{12}=\ldots=c_{1,k+1}=0$ we
see that $$rankC[\{1 \ldots n\},\{2 \ldots k+1\}]=rankC[\{2 \ldots n\},\{2
\ldots
k+1\}],$$ but $C[\{2 \ldots n\},\{2 \ldots k+1\}]=B[\{1 \ldots
n-1\},\{1 \ldots k\}]$
and therefore
\begin{eqnarray}
rankB[\{1 \ldots n-1\},\{1 \ldots k\}]=rankC[\{1 \ldots n\},\{1 \ldots
k+1\}]-1 \label{zero1}
\end{eqnarray}
If $k+1\geq i$ similar reasoning gives
\begin{eqnarray}
&rankB[\{1 \ldots k\},\{1 \ldots (n-1)\}]=rankC[\{1 \ldots (k+1)\},\{1
\ldots n\}]-1 \ \ \ \\
&rankB[\{1 \ldots  k\}]=rankC[\{1 \ldots (k+1) \}]-1
\end{eqnarray}
In particular, $$rankB=rankC-1=rankA-1$$
But if $k+1<i$ we have
\begin{eqnarray}
&rankB[\{1 \ldots k\},\{1 \ldots (n-1)\}]=rankC[\{1 \ldots (k+1)\},\{1
\ldots n\}] \\
&rankB[\{1 \ldots  k\}]=rankC[\{1 \ldots (k+1) \}] \label{zero3}
\end{eqnarray}
Also since $C$ satisfies conditions (\ref{main}) we have
\begin{eqnarray}
rankC[\{1 \ldots (k+1) \}]+(k+1)&\geq&rankC[\{1 \ldots (k+1)\},\{1
\ldots n\}]+ \nonumber \\
 & &rankC[\{1 \ldots n\},\{1 \ldots (k+1)\}]. \label{zero4}
\end{eqnarray}
From (\ref{zero1})--(\ref{zero4}) we conclude that
\begin{eqnarray}
rankB[\{1 \ldots k \}]+k&\geq&rankB[\{1 \ldots k\},\{1 \ldots
(n-1)\}]+\nonumber \\
 & &rankB[\{1 \ldots (n-1)\},\{1 \ldots k\}]. \nonumber
\end{eqnarray} In other words $B$ satisfies conditions (\ref{main}).

We now consider the case when the first column of $A$ is zero, but the
first row is not zero. We let $$C=A^{T}$$ Notice that $C$ satisfies conditions
(\ref{main}) and also the first row of $C$ is zero, but the first
column of $C$ is not zero. It was shown above that there are scalars
$m_{11}$, and $v_{11}$ and a column vector $M_{21}$ and a row vector
$V_{12}$ such that
\begin{eqnarray}
c_{11}=a_{11}&=&m_{11}v_{11} \nonumber \\
C_{21}=A_{12}^{T}&=&M_{21}v_{11} \nonumber \\
C_{12}=A_{21}^{T}&=&m_{11}V_{12} \nonumber
\end{eqnarray}
and additionally that matrix $$C_{22}-M_{21}V_{12}$$ satisfies
conditions (\ref{main}) and $$rank(C_{22}-M_{21}V_{12})=rankC-1=rankA-1$$
We now let
\begin{eqnarray}
l_{11}=v_{11} \nonumber \\
L_{21}=V_{12}^{T} \nonumber \\
u_{11}=m_{11}\nonumber \\
U_{12}=M_{21}^{T}\nonumber
\end{eqnarray}
Equations (\ref{eq:block1})--(\ref{eq:block3}) are obviously satisfied by
this choice of the first column of $L$ and the first row of
$U$. Additionally a factor complement of $A$
$$B=A_{22}-L_{21}U_{12}=(C_{22}-M_{21}V_{12})^{T}$$ satisfies
conditions (\ref{main}) and also $$rankB=rank(C_{22}-M_{21}V_{12})=rankA-1$$

The last case we have to consider is when both the first row and the
first column of $A$ are equal to zero. If $A$ is the zero matrix we
let the first column of $L$ and the first row of $U$ to be equal
zero. Equations (\ref{eq:block1})--(\ref{eq:block3}) are obviously satisfied by
this choice of the first column of $L$ and the first row of $U$. Also
matrix $$B=A_{22}-L_{21}U_{12}$$ is the zero matrix and therefore
obviously satisfies conditions (\ref{main}). In this case $$rankB=rankA=0$$

If $A$ is not the zero matrix let $i$ be the smallest integer such
that the $i$th column of $A$ is not zero. Let us consider the matrix
$C$ such that the first column of $C$ is equal to the $i$th column of
$A$ and columns 2 through $n$ of $C$ are equal to the corresponding
columns of $A$. We want to prove that $C$ satisfies conditions
(\ref{main}). Fix $k$ ($1\leq k \leq n$). If $k$ is less than $i$ then
because the first row of $C$ is zero we have $$rankC[\{1 \ldots k
\},\{1 \ldots n\}] \leq k-1$$ and because columns 2 through $k$ of $C$ are
zero,
but the first column of $C$ is not zero we have $$rankC[\{1 \ldots n
\},\{1 \ldots k\}]=1$$ We conclude that $$rankC[\{1 \ldots k\}]+k\geq
rankC[\{1 \ldots k
\},\{1 \ldots n\}]+rankC[\{1 \ldots n
\},\{1 \ldots k\}]$$ which means that $C$ satisfies condition (\ref{main})
for index $k$.

If $k$ is greater or equal than $i$, then because columns 2 through
$n$ of $C$ are equal to the corresponding columns of $A$  and the
first column of $C$ is equal to the $i$th column of $A$ and also because the
first column of $A$ is zero we have $$rankC[\{1 \ldots k\}]=rankA[\{1
\ldots k\}]$$ In particular $$rankC=rankA$$ Also  $$rankC[\{1 \ldots k
\},\{1 \ldots n\}]=rankA[\{1 \ldots k
\},\{1 \ldots n\}]$$ and  $$rankC[\{1 \ldots n
\},\{1 \ldots k\}]=rankA[\{1 \ldots n
\},\{1 \ldots k\}]$$ Since $A$ satisfies conditions (\ref{main}) for
index $k$, we conclude that $C$ satisfies conditions (\ref{main}) for
index $k$. We already know that because $C$ satisfies conditions (\ref{main})
and the first
row of $C$ is zero, but the first column of $C$ is not zero there are scalars
$m_{11}$, and $v_{11}$ and a column vector $M_{21}$ and a row vector
$V_{12}$ such that
\begin{eqnarray}
c_{11}=0&=&m_{11}v_{11} \nonumber \\
C_{21}&=&M_{21}v_{11} \nonumber \\
C_{12}=A_{12}&=&m_{11}V_{12} \nonumber
\end{eqnarray}
and additionally that matrix $$C_{22}-M_{21}V_{12}$$ satisfies
conditions (\ref{main}) and $$rank(C_{22}-M_{21}V_{12})=rankC-1=rankA-1$$
We now let
\begin{eqnarray}
l_{11}=m_{11} \nonumber \\
L_{21}=M_{21} \nonumber \\
u_{11}=0\nonumber \\
U_{12}=V_{12}\nonumber
\end{eqnarray}
Equations (\ref{eq:block1})--(\ref{eq:block3}) are obviously satisfied by
this choice of the first column of $L$ and the first row of
$U$. Additionally a factor complement of $A$
$$B=A_{22}-L_{21}U_{12}=C_{22}-M_{21}V_{12}$$ satisfies
conditions (\ref{main}) and also
$$rankB=rank(C_{22}-M_{21}V_{12})=rankA-1$$ Note that this case can
also be treated by considering matrix $C$ such that the first {\em
row} of $C$ is equal to the first nonzero row of $A$ and such that
rows 2 through n of $C$ are equal to the corresponding rows of $A$.
$\Box$

We now prove the following well known result

\begin{corollary}
Let $A$ be an n-by-n invertible matrix then $A$ has an LU
factorization iff all principal leading submatrices of $A$ have full
rank.
\end{corollary}

{\em Proof.} Since $A$ is invertible we must have
$$rankA[\{1 \ldots k\},\{1 \ldots n\}]=rankA[\{1 \ldots n\},\{1 \ldots
k\}]=k$$ for all $k=1,\ldots ,n$. Then by theorem (\ref{lu}) $A$ has
an LU factorization iff $$rankA[\{1 \ldots k \}]=k$$ for $k=1, \ldots
,n$ $\Box$

We now give pseudo-code description of the algorithm derived in
theorem (\ref{main}) for finding $L$
and $U$. \\
\\
{\em Algorithm 1}\\
{\bf
Input n-by-n matrix A\\
for k from 1 to n\\
for i from 1 to n \\
for j from i to n \\
if $a_{i,j}\neq 0$  \\
$k$th column of L:= $j$th column of A \\
$k$th row of U:=$a_{i,j}^{-1}$ multiplied by $i$th row of
A\\
GOTO 10\\
elseif $a_{j,i}\neq 0$\\
$k$th column of L:= $i$th column of A \\
$k$th row of U:=$a_{j,i}^{-1}$ multiplied by $j$th row of
A\\
GOTO 10\\
endif\\
endif\\
next j\\
next i\\
$k$th column of L:= the zero vector \\
$k$th row of U:=the zero vector\\
10 A=A-($k$th column of L multiplied by $k$th row of U)\\
next k}\\

We also give a verbal description of the algorithm:\\

\pagebreak

1) Assign priorities to all possible positions in n-by-n matrix. \\

For all k from 1 to n repeat steps 2 through 4\\

2) Find a nonzero element of matrix $A$ which has the smallest
possible integer assigned to its position. If there are no nonzero elements
(i.e if $A$ is the zero matrix) then let $k$th column of $L$  and $k$th row of
$U$
equal to zero and skip to step 4.

3) Suppose the element chosen in step 2 is in position (i,j) then let
$k$th column of $L$ equal to $j$th column of $A$ and $k$th row of $U$
equal to $i$th row of $A$ divided by $a_{i,j}$.

4) Let $$A=A-L[\{1 \ldots n\},\{k\}]U[\{k\},\{1 \ldots n\}]$$

How to assign priorities? \\
We give a simple pseudo-code description and a 4-by-4 example\\
{\bf
let counter=0\\
for i from 1 to n \\
for j from i to n\\
counter=counter+1\\
assign value of counter to (i,j) and (j,i) positions\\
next j\\
next i\\

}

Here is how we would assign priorities to positions in 4-by-4 matrix
$$\left[\begin{array}{cccc}
1&2&3&4\\
2&5&6&7\\
3&6&8&9\\
4&7&9&10
\end{array}\right]$$

We may now state
\begin{theorem}
$A\in M_{n}(F)$ has an LU factorization iff algorithm (1) produces LU
factorization of $A$.
\end{theorem}
{\em Proof.} {\em Necessity.} If $A$ has an LU factorization then by
theorem (\ref{lu}) it satisfies conditions (\ref{main}). Then the
algorithm will produce an LU factorization of $A$. Sufficiency is
trivial. $\Box$

\begin{note}
\label{n1}
\end{note}
Suppose that an n-by-n matrix $A$ satisfies conditions (\ref{main}) and
additionally the first $m$ rows of $A$ are equal to zero and the first
$m$ columns of $A$ are zero. If we apply the algorithm to obtain an LU
factorization of $A$, then the first $m$ rows of $L$ will be equal to
zero and the first $m$ columns of $U$ will be zero.
\begin{note}
\label{n2}
\end{note}
It was shown in the proof of the theorem (\ref{main}) that when we
perform one step of the algorithm (i.e. go from k=p to k=p+1) the rank of
modified  $A$ is less by one than the rank of $A$ before modification,
unless
the rank of $A$ before modification  was already zero.
So if we start with an n-by-n matrix $A$ satisfying conditions (\ref{main})
such that rank deficiency of $A$ is  $m$, then after  (n-m)
steps modified $A$ will be the zero matrix. Therefore the last $m$ columns of
$L$ and the last $m$ row of $U$ will be equal to zero (which means that
the algorithm produces ``full rank'' factorization of $A$). In
particular, if  $A$ has $m$ rows of zeros or $m$ columns of
zeros, then the last $m$ columns of $L$ and rows of $U$ are zero.

\begin{note}
\end{note}

By previous note the ranks of $L$ and $U$ are less or equal to
(n-m). But they also must be greater or equal than the rank of $A$
which is equal to (n-m). Therefore $rankL=rankU=n-m$.

\begin{definition}

An n-by-n matrix $K$ is said to be an almost lower triangular with $m$
extra diagonals if $a(i,j)=0$ whenever $j>i+m$. An almost upper
triangular matrix with $m$ extra diagonals is defined similarly.

\end{definition}

\begin{definition}
We say that an n-by-n matrix $A$ fails conditions(\ref{main}) by at
most (no more than) $m$ if
$$rankA[\{1 \ldots k\}]+k+m\geq rankA[\{1 \ldots k\},\{1 \ldots n\}]+rankA[\{1
\ldots n\},\{1 \ldots k\}]$$

for all $k=1, \ldots,n$.

\end{definition}

Given an n-by-n matrix $A$ we obtain (n+m)-by-(n+m)
matrix $C$ by bordering $A$ with $m$ columns of zeros on the left and
$m$ rows of zeros on the top.

We state a trivial but useful lemma

\begin{lemma}
If $A$ and $C$ are as above then $C$ satisfies conditions (\ref{main})
iff $A$ fails conditions (\ref{main}) by no more than $m$.

\end{lemma}

We now prove
\begin{theorem}
\label{super}
$A\in M_{n}(F)$ can be written as $$A=KW$$ where $K$ is an almost
lower triangular matrix with $m$ extra diagonals and $W$ is an almost
upper triangular matrix with $m$ extra diagonals iff $A$ fails
conditions (\ref{main}) at most by $m$.
\end{theorem}
{\em Proof.} Sufficiency. Given an n-by-n matrix $A$ that fails
conditions (\ref{main}) by no more than $m$ we obtain (n+m)-by-(n+m)
matrix $C$ by bordering $A$ with $m$ columns of zeros on the left and
$m$ rows of zeros on the top. Matrix $C$ satisfies conditions (\ref{main}). We
apply the algorithm to obtain LU
factorization of $C$. As was noted above the first $m$ rows and the
last $m$ columns of $L$ are zero. Also the first $m$ columns and the
last $m$ rows of $U$ are zero. We let $$K=L[\{m+1 \ldots n+m\},\{1
\ldots n\}]$$ and $$W=U[\{1 \ldots n\},\{m+1
\ldots n+m\}]$$ Matrices $K$ and $W$ defined this way have the desired form.
Also block multiplication shows that $$A=KW$$

Necessity. Suppose A can be written as $$A=KW$$ with $K$ and $W$ such
as above. We define $C$, $L$ and $U$ as above. Block multiplication
shows that $$C=LU$$ Therefore by theorem (\ref{lu}) $C$ satisfies
conditions (\ref{main}). It follows that $A$ fails conditions
(\ref{main}) by no more than $m$. $\Box$

\pagebreak

\begin{note}
\label{fr}
\end{note}
If matrix $A$ has rank deficiency $p$ then matrix $C$ defined as above
has rank deficiency $m+p$ and by note (\ref{n2}) the last $m+p$ columns of $L$
and rows of $U$ are equal to zero. Therefore the last $p$ columns of
$K$ are zero and the last $p$ rows of $W$ are zero. In particular our
algorithm produces a full rank factorization of $A$.

\begin{note}
\label{kw}
\end{note}

Notice that algorithm (1) when applied to matrix $C$ such as above
simply ``ignores'' first $m$ zero columns and first $m$ zero
rows. Therefore matrices $K$ and $W$ can be obtained by applying
algorithm (1) to matrix $A$ directly and then letting
$$K=L$$ and $$W=U$$

We now state
\begin{theorem}
An n-by-n matrix $A$ can be written as $$A=KW$$ where $K$ is an almost
lower triangular matrix with $m$ extra diagonals and $W$ is an almost
upper triangular matrix with $m$ extra diagonals iff the algorithm
described in note (\ref{kw}) produces such a factorization.

\end{theorem}
{\em Proof.} Suppose $A$ has such a factorization then by theorem
(\ref{super}) it fails conditions (\ref{main}) by no more than
$m$. Then the algorithm produces a factorization of desired form. Sufficiency
is obvious.$\Box$

\begin{definition}

An n-by-(n+m) matrix $H$ is said to be an almost lower triangular with
$m$ extra columns if the matrix $H[\{1 \ldots n\},\{m+1 \ldots n+m\}]$
is lower triangular. An (n+m)-by-n matrix $W$ is said to be an almost
upper triangular with $m$ extra rows if the matrix $W[\{m+1 \ldots
n+m\},\{1 \ldots n\}]$ is upper triangular.

\end{definition}

\begin{theorem}

\label{long}

An n-by-n matrix $A$ can be written as product of an almost lower
triangular matrix $H$ with $m$ extra columns and an almost upper
triangular matrix $V$ with $m$ extra rows iff $A$ fails conditions
(\ref{main}) by no more than $m$.
\end{theorem}
{\em Proof.}

{\em Necessity.} Suppose $A$ can be written as $$A=HV$$ with $H$ and
$V$ such as above. We obtain an (n+m)-by-(n+m) matrix $C$ by bordering
$A$ with $m$ rows of zeros on the top and $m$ columns of zeros on the
left. We obtain a lower triangular matrix $L$ by bordering $H$ with
$m$ rows of zeros on the top, and an upper triangular matrix $U$ by
bordering $V$ with $m$ columns of zeros on the left. Block
multiplication shows that $$C=LU$$ Therefore by theorem (\ref{lu}) $C$
satisfies conditions (\ref{main}). It follows that $A$ fails
conditions (\ref{main}) by no more than $m$.

{\em Sufficiency.} Suppose $A$ fails conditions (\ref{main}) by no
more than $m$. We define matrix $C$ as above. Matrix $C$ satisfies
conditions (\ref{main}). We now apply algorithm (1) to obtain LU
factorization of $A$. By  note (\ref{n1}) the first $m$ rows of
$L$ are zeros and the first $m$ columns of $U$ are zeros. Block
multiplication shows that $A$ can be written as a product of two
matrices of the desired form $$A=HV$$ with $$H=L[\{m+1 \ldots  n+m \},\{1
\ldots n+m\}]$$ and $$V=U[\{1 \ldots  n+m \},\{m+1
\ldots n+m\}]$$ $\Box$
\begin{note}
\end{note}
Additionally by  note (\ref{n2}) last $m$ columns
of $H$ are zero and last $m$ rows of $V$ are zero.

\begin{note}
\end{note}
Note \ref{fr} shows that factorization of $A$ obtained in this way
is a full rank factorization.

\subsection*{Notes about  computation of LU
factorization in floating point  arithmetic .}

Every invertible matrix which has an LU factorization has a neighborhood in
which every matrix is invertible and has LU factorization. Moreover,
if we require that $u_{ii}=1$ then such a factorization is
unique. Therefore we can define function $A\mapsto [L,U]$. Such a function can
be defined in some neighborhood of any invertible matrix that  has
an LU factorization. The function will be invertible and continuous
in the neighborhood. The inverse will also be 1-1 and
continuous. These facts make it possible to compute LU factorization
of an invertible matrix satisfying conditions (1) using floating point
arithmetic.

 However in general a matrix that has an LU factorization
does not have a neighborhood in which every matrix has an LU
factorization. Also if a factorization exists it does not have to be
unique. Further LU factorization  does not have to depend
continuously on the entries of $A$. Thus, because of the possibility
of the rounding error it is not generally possible  to
compute an LU factorization in the general case in floating point arithmetic
using algorithm (1).

\subsection*{Some Applications.}

\begin{theorem}
Any n-by-n matrix $A$ can be written as $$A=U_{1}LU_{2}$$ with $U_{1}$
and $U_{2}$ upper triangular matrices and $L$ lower triangular matrix.
\end{theorem}
{\em Proof.} By a series of type-3 elementary row operations matrix
$A$ can be transformed into a matrix $C$ such that for all $k=1,
\ldots ,n$ we have $$rankC[\{1 \ldots k\}]=rankC[\{1 \ldots n\},\{1
\ldots k\}]$$ In particular, since $$rankC[\{1 \ldots k\},\{1 \ldots
n\}]\leq k$$ we see that $C$ satisfies conditions (\ref{main}).
 In fact one can transform $A$ into $C$ using only elementary type-3 row
operations that add
multiples of $i$th row to the $j$th row  with $j<i$. Any
series of such operations can be realized by multiplying $A$ by an
invertible upper triangular matrix $U$ on the left. We have $$C=UA$$
but by theorem (\ref{lu}) $$C=LU_{2}$$ with $L$ lower triangular and
$U_{2}$ upper triangular. Since $U$ is invertible we can let
$$U_{1}=U^{-1}$$  That gives us $$U_{1}UA=A=U_{1}LU_{2}$$ $\Box$
\begin{note}
\end{note}
Notice that $U_{1}$ can be taken invertible. Also because $$rankC[\{1 \ldots
k\}]=rankC[\{1 \ldots n\},\{1
\ldots k\}]$$  it
follows  that $L$ can be taken invertible [ see LM ].

\begin{corollary}
Any n-by-n matrix $A$ can be written as $$A=L_{1}UL_{2}$$ with $L_{1}$
and $L_{2}$ lower triangular matrices and $U$ an upper triangular matrix.

\end{corollary}

{\em Proof.} It follows from the fact that $A^{T}$ can be written as
$$A^{T}=U_{1}LU_{2}$$ with $U_{1}$
and $U_{2}$ upper triangular matrices and $L$ lower triangular
matrix. $\Box$

\begin{note}
\end{note}

Notice that $U$ and $L_{2}$ can be taken invertible.

We prove the following well known result
\begin{theorem}
Any n-by-n matrix $A$ can be written as $$A=PLU$$ with $U$
an upper triangular matrix, $L$ a lower triangular matrix and $P$
permutation matrix.
\end{theorem}

{\em Proof.} Any matrix $A$ can be multiplied on the left by such
permutation matrix $P_{0}$ that matrix $$C=P_{0}A$$ satisfies the
following equality for all $k=1, \ldots ,n$ $$rankC[\{1 \ldots k\}]=rankC[\{1
\ldots n\},\{1
\ldots k\}]$$ In particular, since $$rankC[\{1 \ldots k\},\{1 \ldots
n\}]\leq k$$ we see that $C$ satisfies conditions
(\ref{main}). therefore by theorem (\ref{main}) $C$ can be written as
$$C=LU$$ with a lower triangular $L$ and an upper triangular $U$. We
let $$P=P_{0}^{-1}$$ That gives us $$A=PP_{0}A=PC=PLU$$ $\Box$.

\begin{corollary}
Any n-by-n matrix $A$ can be written as $$A=LUP$$ with $U$
an upper triangular matrix, $L$ a lower triangular matrix and $P$
permutation matrix.
\end{corollary}
{\em Proof.} It follows from the fact that $$A^{T}=P_{0}L_{0}U_{0}$$
with $U_{0}$
an upper triangular matrix, $L_{0}$ a lower triangular matrix and $P_{0}$
permutation matrix. Indeed, we can write $$A=LUP$$ with
$$L=U_{0}^{T}$$ $$U=L_{0}^{T}$$ and $$P=P_{0}^{T}$$ $\Box$

\end{document}